\documentclass[a4paper]{article}
\usepackage[dvips]{graphicx}
\usepackage{amssymb}
\usepackage{amsmath}
\usepackage{caption}
\usepackage{hyperref}

\makeatletter
\newcommand*{\centerfloat}{%
  \parindent \z@
  \leftskip \z@ \@plus 1fil \@minus \textwidth
  \rightskip\leftskip
  \parfillskip \z@skip}
\makeatother

\newcommand {\aplt} {\ {\raise-.5ex\hbox{$\buildrel<\over\sim$}}\ }
\makeatletter
\def\dddots{\mathinner{\mkern1mu\raise\p@
    \hbox{.}\mkern2mu\raise4\p@\hbox{.}\mkern2mu
    \raise7\p@\vbox{\kern7\p@\hbox{.}}\mkern1mu}}%
\makeatother

\newtheorem{Proposition}{Proposition}

\newtheorem{Remark}{Remark}

\newenvironment{AMS}{\small\bf 2010 AMS subject classification: }{}

\newfont{\BBB}{msbm10 scaled \magstep1}
\newfont{\BBS}{msbm10}

\begin{document} 

\title{\bf Discrete norming inequalities   
on sections of sphere, ball and torus\thanks{Work partially 
supported by the BIRD163015, PRAT-CPDA143275 and DOR funds of 
the University of
Padova, and by the GNCS-INdAM. This research has been accomplished 
within the RITA
``Research ITalian network on Approximation''.}}

\author{Alvise Sommariva$^1$ and Marco Vianello$^1$}

\maketitle

\footnotetext[1]{Department of Mathematics, University of
Padova, Italy\\corresponding author: marcov@math.unipd.it}

\begin{abstract} 
By discrete trigonometric norming inequalities on 
subintervals of 
the period, we construct   
norming meshes with optimal cardinality growth for algebraic 
polynomials on sections of sphere, ball and torus. 
 
\end{abstract}
\vskip0.2cm
\noindent
\begin{AMS}
{\rm 41A10, 41A63, 42A05, 65D18.}
\end{AMS}

\section{Introduction}
Polynomial inequalities based on the notion of polynomial 
mesh have been recently playing a relevant role in multivariate 
approximation theory, as well in its computational applications.   

We recall that a {\em polynomial mesh} of a  
compact subset $K$ of a manifold 
$\mathcal{M}\subseteq \mathbb{R}^d$,      
is a sequence of finite norming subsets $\mathcal{A}_n\subset K$ such that 
the polynomial inequality
\begin{equation} \label{polmesh}
\|p\|_K\leq 
c\,\|p\|_{\mathcal{A}_n}\;,\;\forall
p\in \mathbb{P}_n^d(K)\;,
\end{equation}
holds for some $c>1$ independent of $p$ and $n$,  
where 
$card(\mathcal{A}_n)=\mathcal{O}(N^s)$, $N=dim(\mathbb{P}_n^d(K))$, 
$s\geq 1$. 

Here and below we denote by $\mathbb{P}_n^d(K)$ the subspace
of $d$-variate real polynomials of total degree not exceeding $n$
restricted to $K$, and by
$\|f\|_X$ the sup-norm of a bounded real function on
a discrete or continuous compact set $X\subset \mathbb{R}^d$.

Following \cite{CL08}, when $c=c_n$ depends on $n$ but with 
subexponential growth, we speak of a {\em weakly admissible} polynomial 
mesh. In \cite{CL08}, in the case of $c$ independent of $n$ the polynomial 
mesh is termed {\em admissible}. In this paper we focus on admissible 
polynomial meshes, that we simply term ``polynomial meshes'' as in 
(\ref{polmesh}). 

Observe that $\mathcal{A}_n$ is 
$\mathbb{P}_n^d(K)$-determining (i.e., a
polynomial vanishing there 
vanishes
everywhere on $K$), 
consequently $card(\mathcal{A}_n)\geq N$. A polynomial mesh may  
then be termed {\em optimal} when 
$s=1$.  

All these notions can be given 
more generally for $K\subset \mathbb{C}^d$ but we restrict 
here to real compact sets. They are extensions of notions usually given 
for $\mathcal{M}=\mathbb{R}^d$. 

Polynomial meshes were formally introduced in the seminal paper 
\cite{CL08} and then studied from both the theoretical 
and the computational point of view. Among their  
features, we recall that they:  
\begin{itemize}   
\item are affinely invariant and can be extended 
by algebraic transforms, finite union and product 
\cite{BCLSV11,CBC16,CL08};
\item are stable under small 
perturbations \cite{PV13}; 
\item contain near optimal
interpolation sets of Fekete type (maximal Vandermonde determinant), 
computable by standard numerical linear algebra algorithms  
\cite{BCLSV11,BDMSV10,CL08}; 
\item are near optimal for uniform Least Squares 
approximation \cite{BCLSV11,CL08}, namely
\begin{equation} \label{LSnorm}
\|\mathcal{L}_{\mathcal{A}_n}\|
=\sup_{f\in C(K),f\neq 0}\frac{\|\mathcal{L}_{\mathcal{A}_n}f\|_K}
{\|f\|_K}\leq c\,\sqrt{card(\mathcal{A}_n)}\;,
\end{equation}
where $\mathcal{L}_{\mathcal{A}_n}$ is the $\ell^2(\mathcal{A}_n)$-orthogonal
projection operator
$C(K)\to \mathbb{P}_n^d(K)$ (the discrete LS operator at $\mathcal{A}_n$),
from which easily follows
\begin{equation} \label{LSerr}
\|f-\mathcal{L}_{\mathcal{A}_n}f\|_K\leq \left(1+c\,
\sqrt{card(\mathcal{A}_n)}\right)\,
\min_{p\in \mathbb{P}_{n}^d(K)}\|f-p\|_K\;;
\end{equation}
 
\item can be used in pluripotential numerics, e.g. for computing the 
pluripotential Green function 
and the multivariate transfinite diameter \cite{P17}; 
\item can be used in the framework of polynomial optimization 
\cite{PV17-1,V17}.  
\end{itemize}

Optimal polynomial meshes have been constructed on several classes of compact 
sets, such as polygons and polyhedra, circular and spherical 
sections, convex bodies and 
star-shaped domains, general compact domains with regular boundary, by 
different 
analytical and geometrical techniques; we refer 
the reader, e.g., to  
\cite{BV12,CBC16,CL08,GSV16,JSW99,K11,K17,LSV16,P16,SV15} and 
the 
references 
therein, for a comprehensive view of construction methods and 
applications. 

\section{Norming inequalities on surface/solid sections}
In this paper we survey several cases where $K$ is the image 
of a box of scalar and angular variables, by a 
geometric transformation whose components are in tensor product 
spaces of algebraic and trigonometric polynomials of degree one. 
The angular variables are divided into periodic ones (the relevant 
intervals have length $2\pi$) and superiodic ones (the relevant
intervals have length $<2\pi$). In the sequel, we denote by 
\begin{equation} \label{trig}
\mathbb{T}_n([u,v])=span(1,cos(\theta),\sin(\theta),\dots,
\cos(n\theta),\sin(n\theta))\;,\;\;\theta\in [u,v]\;,
\end{equation}
the space of univariate 
trigonometric 
polynomials 
of degree not exceeding $n$, restricted to the angular interval 
$[u,v]$. When $v-u<2\pi$ we are in a {\em subperiodic}  
instance (a subinterval of the period). It is worth recalling that 
trigonometric approximation on subintervals of the period, also termed 
theory of ``Fourier extensions'' in some contexts, 
has been object 
of several studies in the recent literature, cf. e.g. 
\cite{AHV14,BV12-2,DFV12} 
with the references therein.  

More precisely, we consider compact sets of the form   
$$
K=\sigma(I\times \Theta)\;,\;\;
\sigma=\left(\sigma_\ell(t,\theta)\right)_{1\leq \ell \leq d}\;,
$$
\begin{equation} \label{sigma}
t\in I=I_1\times\dots\times I_{d_1}\;,\;\;
\theta\in \Theta=\Theta_1\times\dots\times 
\Theta_{d_2+d_3}\;,
\end{equation}
\begin{equation} \label{prodtens}
\sigma_\ell\in
\bigotimes_{i=1}^{d_1}\mathbb{P}_1(I_i)
\otimes \bigotimes_{j=1}^{d_2+d_3}\mathbb{T}_1(\Theta_j)
\;\;,\;\;1\leq \ell \leq d\;, 
\end{equation}
where $d_1,d_2,d_3\geq 0$, 
and $I_i=[a_i,b_i]$, $1\leq i\leq d_1$ (algebraic 
variables), $\Theta_j=[u_j,v_j]$ with 
$v_j-u_j=2\pi$, $1\leq 
j\leq d_2$ (periodic trigonometric variables) and $v_j-u_j<2\pi$, 
$d_2+1\leq
j\leq d_2+d_3$ (subperiodic trigonometric variables).

A number of arcwise sections of disk, sphere, ball, surface and 
solid torus, fall into the class (\ref{sigma})-(\ref{prodtens}). 
For example, a {\em circular 
sector} 
of the unit disk with angle $2\omega$, $\omega<\pi$, corresponds 
up to a rotation to  
$d_1=d_3=1$, $d_2=0$, 
\begin{equation} \label{circsect}
\sigma(t,\theta)
=(t\cos(\theta),t\sin(\theta))\;,\;\;(t,\theta)\in [0,1]\times 
[-\omega,\omega]\;,
\end{equation}
(polar coordinates). 
Similarly, 
a {\em circular segment} with angle $2\omega$ 
(one of the two portions of the disk cut by a line) corresponds 
up to a rotation again to 
$d_1=d_3=1$, $d_2=0$, but now
\begin{equation} \label{circsegm} 
\sigma(t,\theta)=(\cos(\theta),t\sin(\theta))\;,\;\;
(t,\theta)\in [-1,1]\times [-\omega,\omega]\;. 
\end{equation}
On the other hand, a {\em rectangular tile} of a torus 
corresponds in our notation to
$d_3=2$, $d_1=d_2=0$,
\begin{equation} \label{torustile}
\sigma(\theta)=\left((R+r\cos(\theta_1))\cos(\theta_2),(R+r\cos(\theta_1))
\sin(\theta_2),r\sin(\theta_1)\right)\;,
\end{equation}
$\theta=(\theta_1,\theta_2)\in [\omega_1,\omega_2]\times 
[\omega_3,\omega_4]$,  
where $R$ and $r$ are  
the major and
minor radius of the torus. In the degenerate case $R=0$ we get a 
so-called {\em 
geographic rectangle} 
of a sphere of radius $r$, i.e. the region comprised between 
two given 
latitudes and longitudes.
We refer the reader to 
\cite{DFSV13,DFV14,GSV16,SV15} for several planar and surface 
sections of this kind. Examples of solid sections will be given 
below.     

The key observation in order to use the geometric structure 
to construct polynomial meshes on compact sets in the class 
(\ref{sigma})-(\ref{prodtens}), is that if $p\in \mathbb{P}_n^d(K)$ 
then 
\begin{equation} \label{prodtens2}
p\circ \sigma\in
\bigotimes_{i=1}^{d_1}\mathbb{P}_n(I_i)
\otimes \bigotimes_{j=1}^{d_2+d_3}\mathbb{T}_n(\Theta_j)\;.
\end{equation}
Indeed, in the univariate case Chebyshev-like optimal polynomial meshes 
are known 
for both algebraic polynomials and trigonometric polynomials on 
intervals (even 
on subintervals of the period). This allows to prove the 
following: 
\begin{Proposition}
Let $K\subset \mathbb{R}^d$ be a compact set of the form 
(\ref{sigma})-(\ref{prodtens}). Then, for every fixed $m>1$, 
$K$ possesses a polynomial mesh 
$\mathcal{A}_n=\mathcal{A}_n(m)$ (cf. (\ref{polmesh})) such that
\begin{equation} \label{tensmesh}
c=\alpha^{d_1+d_2}\,\beta^{d_3}\;,\;\;card(\mathcal{A}_n)\leq 
N_1^{d_1}\,N_2^{d_2+d_3}\;,
\end{equation}
where
\begin{equation} \label{c1c2}
\alpha=\alpha(m)=\frac{1}{\cos(\pi/(2m))}\;,\;\;
\beta=\beta(m)=\frac{m}{m-1}\;,
\end{equation}
\begin{equation} \label{cards}
N_1=\lceil mn+1\rceil\;,\;\;N_2=\lceil 2mn+1\rceil\;.
\end{equation}

\end{Proposition}
\vskip0.3cm
\noindent
{\bf Proof.\/} 
In view of the tensorial structure in (\ref{prodtens2}), 
we restrict our attention to univariate instances. 
Consider 
the Chebyshev-Lobatto nodes
of an interval $[a,b]$ (via the affine transformation
$\mathcal{A}(s)=\frac{b-a}{2}\,s+\frac{b+a}{2}$),
\begin{equation} \label{cheblob}
X_\nu([a,b])=\left\{\mathcal{A}(\xi_j)\right\}
\subset [a,b]\;,\;\;
\xi_j=\cos{(j\pi/\nu)}\;,\;0\leq j \leq \nu\;,
\end{equation}
the classical Chebyshev nodes (the zeros of $T_{\nu+1}(s)$) 
\begin{equation} \label{chebzeros}
Z_\nu([a,b])=\left\{\mathcal{A}(\eta_j)\right\}\subset
(a,b)\;,\;\;\eta_j=\cos\left(\frac{(2j+1)\pi}{2(\nu+1)}\right)\;,\;\;0\leq j
\leq \nu\;,
\end{equation}
and the Chebyshev-like ``subperiodic'' angular nodes
\begin{equation} \label{chebangles}
W_\nu([u,v])=\psi_\omega(Z_{2\nu}([-1,1]))+\frac{u+v}{2}
\subset (u,v)\;,\;\;\omega=\frac{v-u}{2}\leq \pi\;,
\end{equation}
obtained by the nonlinear tranformation
$\psi_\omega(s)=2\arcsin\left(\sin\left(\frac{\omega}{2}\right)
\,s\right)$, $s\in [-1,1]$.
Notice that $card(X_\nu)=card(Z_\nu)=\nu+1$, $card(W_\nu)=2\nu+1$. 
Moreover, all the nodal families cluster at the interval endpoints, except  
for the periodic case $W_\nu([u,u+2\pi])$, where the nodes are equally 
spaced. 

These nodal sets satisfy the following fundamental 
norming inequalities 
$$
\|g\|_{[a,b]}\leq \alpha\,\|g\|_{X_{mn}}\;,\;\;\forall g\in 
\mathbb{P}_n([a,b])\;,
$$
$$
\|\tau\|_{[u,v]}\leq \alpha\,\|\tau\|_{W_{mn}}\;,\;\;\forall \tau\in
\mathbb{T}_n([u,v])\;,\;\;v-u=2\pi\;,
$$
\begin{equation} \label{univpol}
\|\tau\|_{[u,v]}\leq \beta\,\|\tau\|_{W_{mn}}\;,\;\;\forall \tau\in
\mathbb{T}_n([u,v])\;,\;\;v-u<2\pi\;,
\end{equation}
where $\alpha$ and $\beta$ are defined in (\ref{c1c2}).
The first
and second inequality are well-known results of polynomial approximation 
theory 
obtained
by Ehlich and Zeller \cite{EZ64}, the
third
has been recently proved in the framework of subperiodic trigonometric
approximation \cite{PV16,V14}, improving previous $\omega$-dependent
estimates in \cite{K11}).

By (\ref{prodtens2}) we then obtain
\begin{equation} \label{tenspolmesh}
\|p\|_K=\|p\circ \sigma\|_{I\times \Theta}
\leq \alpha^{d_1+d_2}\,\beta^{d_3}\|p\circ 
\sigma\|_{\mathcal{B}_n}=\alpha^{d_1+d_2}\,\beta^{d_3}
\|p\|_{\sigma(\mathcal{B}_n)}\;,
\end{equation}
where
\begin{equation} \label{prodmesh}
\mathcal{B}_n=\mathcal{B}_n(m)=\left(X_{mn}(I_1)\times\dots \times 
X_{mn}(I_{d_1})\right)
\times \left(W_{mn}(\Theta_1)\times\dots \times 
W_{mn}(\Theta_{d_2+d_3})\right)\;.
\end{equation}
Observe that $card(\mathcal{B}_n)=\lceil mn+1\rceil ^{d_1}\lceil 2mn+1\rceil 
^{d_2+d_3}$. 
The cardinality estimate in (\ref{tensmesh}) then follows with 
$\mathcal{A}_n=\sigma(\mathcal{B}_n)$, since $\sigma$ in general is not 
injective. Notice, finally, that $X_{mn}$ could be substituted by 
$Z_{mn}$ in all the construction.$\;\;\;\;\square$
\vskip0.5cm    

\begin{Remark}
{\em 
If $N=dim(\mathbb{P}_n^d(K))\sim \gamma n^\lambda$, $n\to \infty$ 
($\lambda\leq d$), 
the polynomial mesh (\ref{tensmesh}) is {\em optimal} when 
$d_1+d_2+d_3=\lambda$, 
since then $card(\mathcal{A}_n)=\mathcal{O}(n^\lambda)=\mathcal{O}(N)$. 
This happens for all the sections of sphere, ball and torus 
listed below.
\/}
\end{Remark}

\begin{Remark}
{\em If besides $n$ also $m$ is an integer, the mesh $\mathcal{A}_n(m)$ 
can be considered in three ways: 
\begin{itemize}
\item as an optimal (admissible) polynomial mesh for degree $n$ with  
$c=c(m)=(\alpha(m))^{d_1+d_2}\,(\beta(m))^{d_3}$ (or symmetrically for 
degree 
$m$ with $c=c(n)$);
\item as a weakly admissible polynomial mesh for degree $mn$ with
$$ 
c_{mn}=\mathcal{O}\left((\log(mn))^{d_1+d_2+d_3}\right)
$$
in view 
of the results in \cite{BV12,DFV13}, see also \cite{SV15}.
\end{itemize}
$\;$
\/} 
\end{Remark}
\vskip0.3cm
\noindent

\subsection{Planar sections}
Several sections of the disk (ellipse) can be described 
in a unifying way as {\em linear blending} of arcs, namely 
by 
\begin{equation} \label{blending}
\sigma(t,\theta)=t\,P_1(\theta)+(1-t)\,P_2(\theta)\;,\;\;t\in 
[0,1]\;,\;\;\theta\in [\omega_1,\omega_2]\;,\;\;\omega_2-\omega_1<2\pi\;,
\end{equation}  
with
\begin{equation} \label{Pi}
P_i(\theta)=A_i\cos(\theta)+B_i\sin(\theta)+C_i\;,\;\;i=1,2\;,
\end{equation}
where $A_i,B_i,C_i$ are suitable 2-dimensional vectors (with $A_i$ and 
$B_i$ not both null). The transformation (\ref{blending}) has the form 
(\ref{prodtens}), with $d=2$ and $d_1=d_3=1$, $d_2=0$. Among such sections 
we may quote circular symmetric or even asymmetric 
{\em 
sectors}, circular {\em annuli}, {\em segments}, {\em zones}, {\em 
lenses}. We recall 
that a zone is the section of a disk cut by two parallel lines, whereas a 
lens is the intersection of two overlapping disks. 
 
We stress that there are in general different representations
of the form (\ref{blending})  
for a given circular section. For example, a  
circular segment (of the unit disk) has 
at least four blending representations, namely 
taking for example a semi-angle $\omega$ 
we may have 
$P_1(\theta)=(cos(\theta),\sin(\theta))$, and  
$P_2(\theta)=(\cos(\theta),\sin(\pi/2-\omega))$
(arc-segment blending see Fig. 1 top-left), or  
$P_2(\theta)=(0,\sin(\pi/2-\omega))$ 
(arc-point blending that is generalized sector, see Fig. 1 top-right), 
both with $\theta\in
[\pi/2-\omega,\pi/2+\omega]$,  
or $P_2(\theta)=(-\cos(\theta),\sin(\theta))$ 
(arc-arc blending) with $\theta\in
[\pi/2-\omega,\pi/2]$ (Fig. 1 bottom-left),  
or with $\theta\in
[\pi/2-\omega,\pi/2+\omega]$ (Fig. 1 bottom-right).  
Notice that with the last choice the blending transformation $\sigma$ is 
not injective and the number of points is essentially halved, thanks to  
the transformation symmetry.  
We refer the reader to
\cite{DFSV13,DFV12,SV15} for these and several other examples of blending
type
in the framework of interpolation and cubature.

\begin{figure}[!ht] \centering
\includegraphics[scale=0.30,clip]{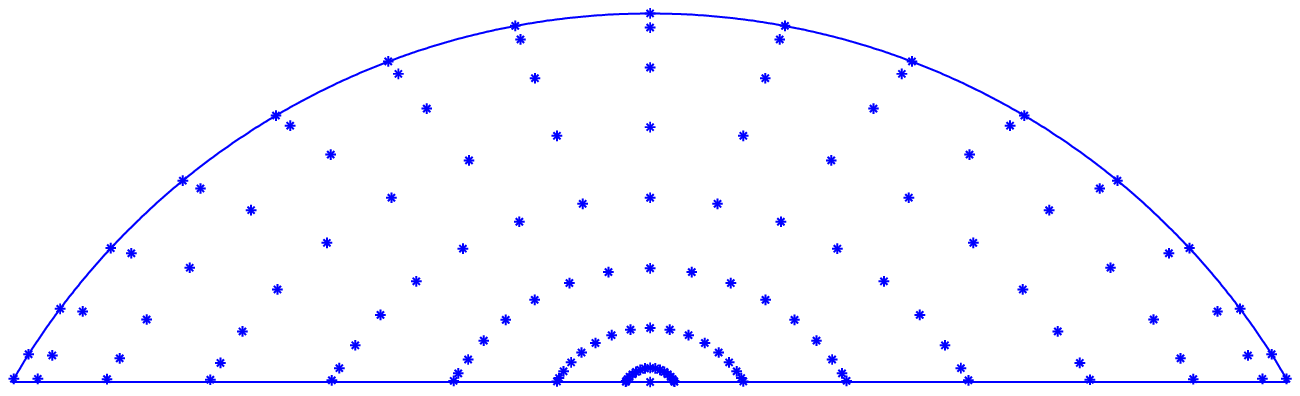}\hfill   
\includegraphics[scale=0.30,clip]{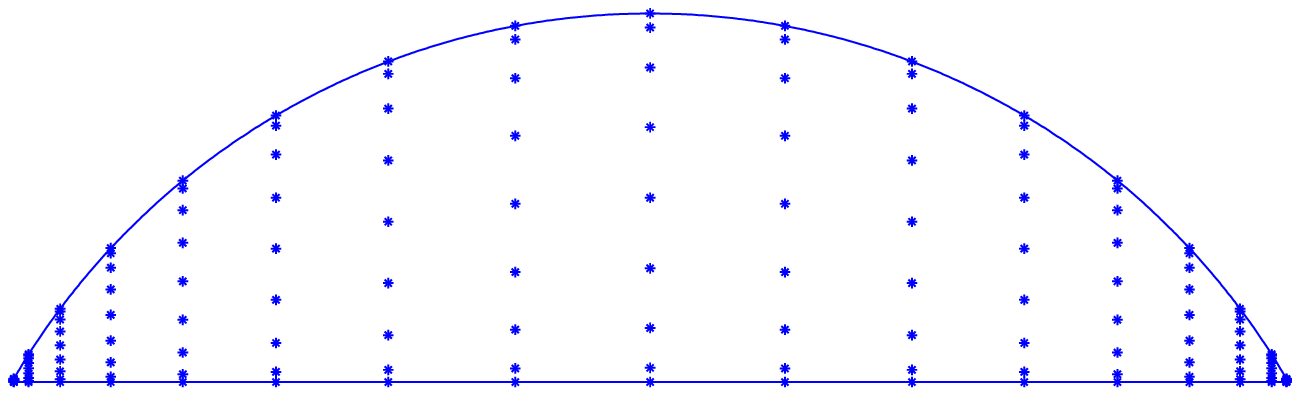}\hfill
\includegraphics[scale=0.30,clip]{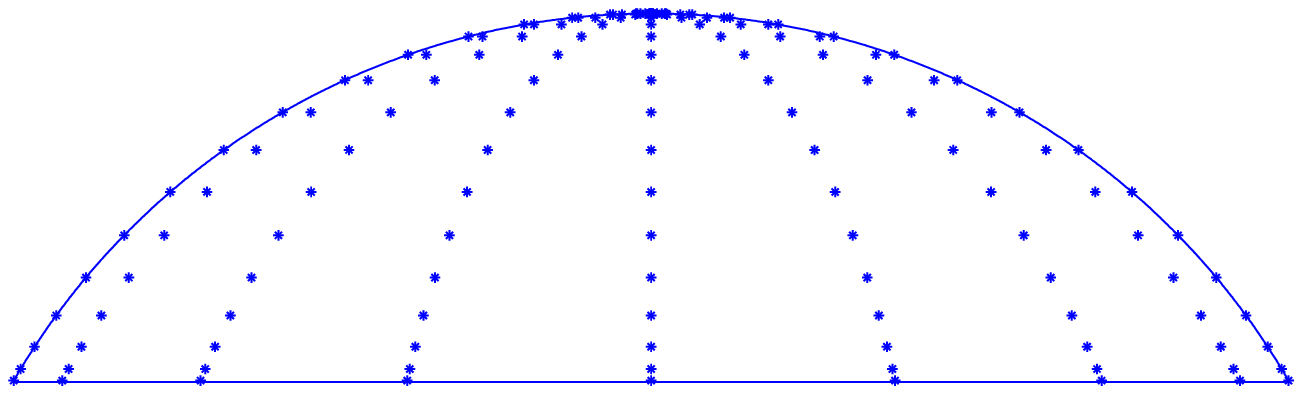}\hfill
\includegraphics[scale=0.30,clip]{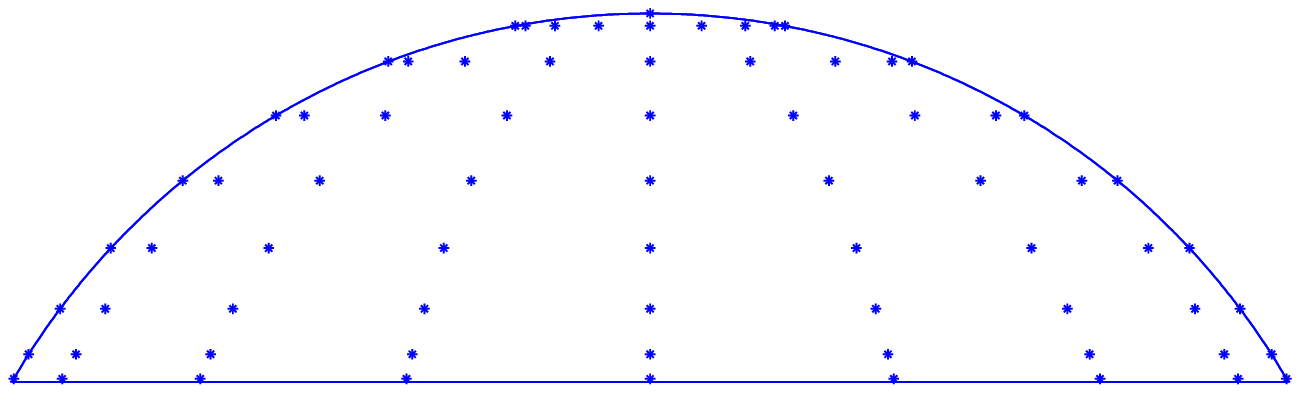}\hfill
\caption{Four blending polynomial meshes for degree $n=4$ on a circular 
segment with angle $2\pi/3$, $m=2$ and $c=2\sqrt{2}$; the mesh 
cardinality is $9\times 17$ 
(top-left, top-right and bottom-left) and $9\times 9$ (bottom-right).}
\label{blendsegm}   
\end{figure}

The 
corresponding polynomial meshes 
\begin{equation} \label{blendmesh}
\mathcal{A}_n(m)=\sigma\left(X_{mn}([0,1])\times
W_{mn}([\omega_1,\omega_2])\right)\;,
\end{equation}
are optimal, since $d_1+d_2+d_3=2$
and
$N=dim(\mathbb{P}_n^2(K))=(n+1)(n+2)/2\sim n^2/2$ (cf. Remark 1), 
with $c=\alpha \beta$.

A special role is played by planar {\em circular lunes} 
(difference of a disk with a second    
overlapping disk), that do not fall in the arc blending class and 
correspond to a transformation $\sigma$ where two angular 
variables are involved. Indeed a lune of the unit disk, whose boundary is 
given by two 
circular 
arcs, a longer one with semiangle say $\omega_2$ and a shorter one 
with semiangle say $\omega_1$, can be described by the bilinear 
trigonometric transformation $\sigma=(\sigma_1,\sigma_2)$   
$$
\sigma_1(\theta_1,\theta_2)=\cos(\theta_2)-\frac{\cos(\omega_1)}
{\sin(\omega_1)}\,\sin(\theta_2)+\frac{1}
{\sin(\omega_1)}\,\cos(\theta_1)\sin(\theta_2)\;,
$$
\begin{equation} \label{lune}
\sigma_2(\theta_1,\theta_2)=\frac{1}
{\sin(\omega_1)}\,\cos(\theta_1)\sin(\theta_2)\;,
\;\;(\theta_1,\theta_2)\in 
[-\omega_1,\omega_1]\times [0,\omega_2]\;.
\end{equation} 
As proved in \cite{DFV14} such a transformation maps (not injectively 
since $\sigma_1(\theta_1,0)\equiv 1$) the boundary of the rectangle onto 
the boundary of the lune
(preserving the orientation) and has positive Jacobian, so
that it is a diffeomorphism of the interior of the rectangle onto the 
interior of the lune.
Observe that (\ref{lune}) fall into the class (\ref{prodtens}) 
with $d=2$, $d_1=d_2=0$, $d_3=2$. 
The corresponding polynomial meshes 
\begin{equation} \label{lunemesh}
\mathcal{A}_n(m)=\sigma\left(W_{mn}([-\omega_1,\omega_1])\times
W_{mn}([0,\omega_2])\right)\;,
\end{equation}
are again optimal, with $c=\beta^2$.

\subsection{Sections of sphere and torus}
The relevant transformation $\sigma$ has the form (\ref{torustile}), 
which 
characterizes suitable arcwise sections of the sphere ($R=0$) or of the 
torus 
($R>0$). In this surface instances we have $d_1=0$ (only angular 
coordinates are involved). As observed above, a {\em spherical or 
toroidal 
rectangle} corresponds to $d_2=0$, $d_3=2$, $\theta\in 
[\omega_1,\omega_2]\times [\omega_3,\omega_4]$, where 
$\omega_2-\omega_1,\omega_4-\omega_3>0$,  
$[\omega_1,\omega_2]\subseteq [-\pi/2,\pi/2]$ and 
$[\omega_3,\omega_4]\subset 
[-\pi,\pi]$. The special case $[\omega_1,\omega_2]=[-\pi/2,\pi/2]$ 
is a so-called {\em spherical lune}.   

Special instances may be periodic in one of the variables. For example, 
a {\em spherical collar} corresponds to
$d_2=d_3=1$, $[\omega_3,\omega_4]=[-\pi,\pi]$ 
(the 
whole longitude range). Geometrically it is the portion of sphere 
between two parallel cutting planes. 

The same holds for a {\em 
spherical 
cap}, one of the two portions of a sphere cut by a plane. Focusing on 
a polar cap (up to a rotation) we have now   
$[\omega_1,\omega_2]=[\frac{\pi}{2}-\omega,\frac{\pi}{2}+\omega]$ with 
$\omega\leq 
\pi/2$ (the 
latter is a generalization of the usual latitude). 

For $R>0$ 
and the same 
angular intervals we obtain what we may call a {\em toroidal collar} 
and a {\em toroidal cap} (geometrically, we are cutting a standard torus 
by planes parallel to the $x_1x_2$-plane). 

On the other hand, for 
$[\omega_3,\omega_4]\subset [-\pi,\pi]$ and 
$[\omega_1,\omega_2]=[-\pi,\pi]$ 
we get a (surface) {\em toroidal slice} (we are cutting the 
torus by two half-planes hinged on the $x_3$-axis).     

All the corresponding polynomial meshes 
\begin{equation} \label{surfmesh}
\mathcal{A}_n(m)=\sigma\left(W_{mn}([\omega_1,\omega_2])\times 
W_{mn}([\omega_3,\omega_4])\right)\;,
\end{equation} 
are optimal, with $c=\beta^2$ (or $c=\alpha\beta$ in periodic 
instances in one of the variables).  
In fact, $d_1+d_2+d_3=2$ 
and 
$N=dim(\mathbb{P}_n^3(K))\sim \gamma n^2$, specifically $N=(n+1)^2$ 
(sphere) or $N=2n^2$ (torus); cf. Remark 1. More generally, if $K$ is a 
polynomial 
determining compact subset of a {\em real algebraic variety} 
$\mathcal{M}\subset \mathbb{R}^d$ defined as the zero set of an 
irreducible real polynomial 
of degree $k$, 
then 
\begin{equation} \label{variety}
N=dim(\mathbb{P}_n^d(K))=dim(\mathbb{P}_n^d(\mathcal{M}))={n+d\choose 
d}-{n-k+d\choose d}\;,
\end{equation}
for $n\geq k$, and thus $N\sim \frac{k}{(d-1)!}\,n^{d-1}$ for $n\to 
\infty$ ($d$ and $k$ fixed). Indeed, the sphere is a quadric 
($k=2$) whereas the torus is a quartic ($k=4$) surface in $d=3$ variables; 
see, e.g., \cite{CLOS15} for the relevant algebraic geometry notions.    

In Figure 2 we show two examples of surface optimal polynomial meshes 
on sections of sphere and torus. The numerical
codes that
generate the optimal polynomial meshes are available at \cite{SV17}.

\begin{figure}[!ht] \centering
\includegraphics[scale=0.36,clip]{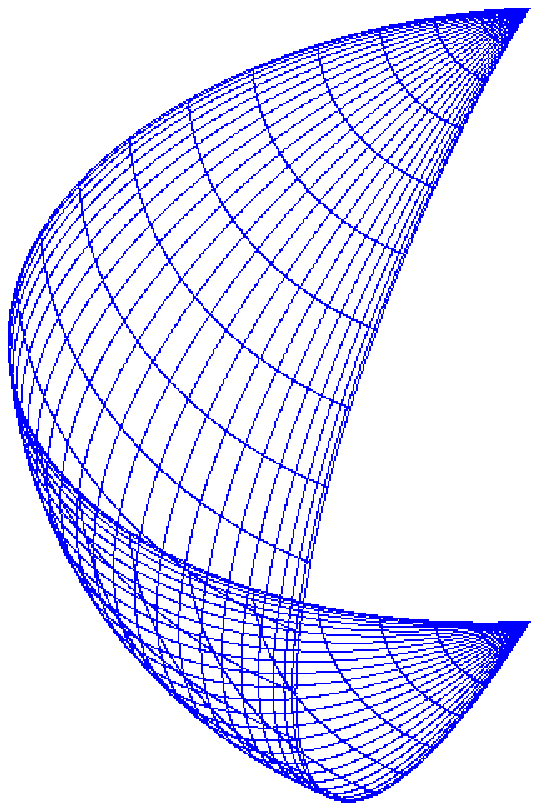}\hfill
\includegraphics[scale=0.36,clip]{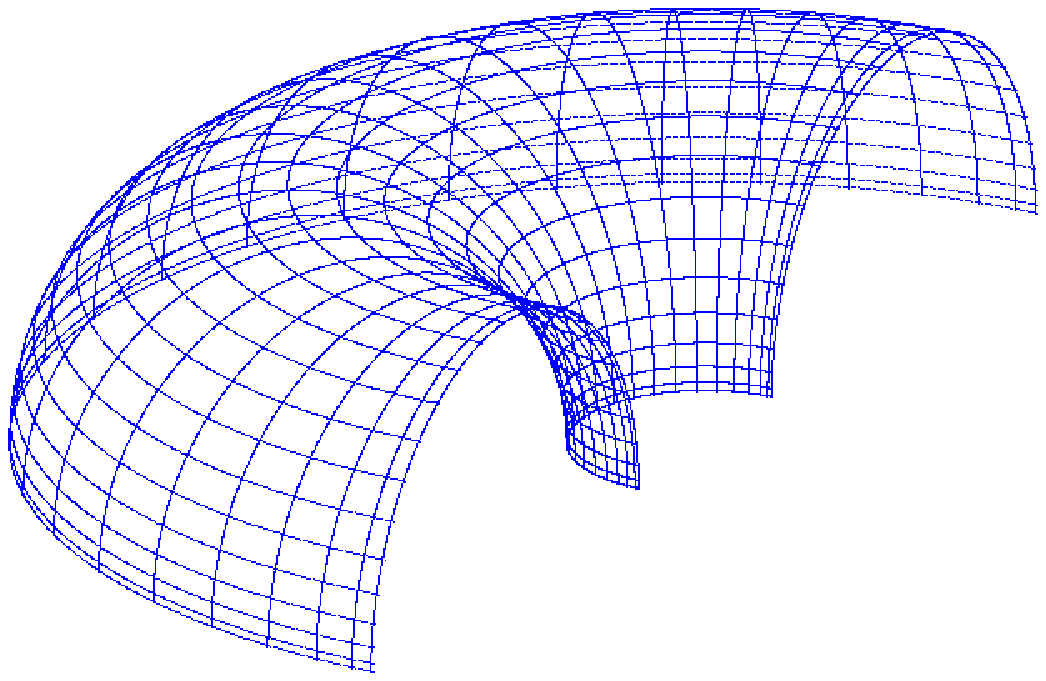}\hfill
\caption{Optimal polynomial meshes (bullets) for $n=4$ with $m=3$  
on a spherical lune (left) and a toroidal cap section 
(right).} 
\label{surfmeshes} \end{figure}

\subsection{Sections of ball and solid torus} 
We focus now on solid arcwise sections of ball and torus. A first 
class of sections corresponds to a rotation of planar sections around 
a coplanar axis, 
by an angle possibly smaller than $2\pi$. Such ``solids of rotation'', 
together 
with the corresponding polynomial meshes, can be conveniently described   
using {\em generalized cylindrical coordinates}
\begin{equation} \label{cyl}
x=(x_1,x_2,x_3)=(r\cos(\phi),r\sin(\phi),z)\;,
\end{equation} 
where we have taken with no loss of generality the $x_3$-axis 
as the rotation axis, we have named for convenience $\phi$ the rotation 
angle, $z$ the $x_3$ coordinate and $r$ assumes also negative values 
(namely $(r,\phi)$ are generalized polar coordinates in planes orthogonal 
to the rotation axis). Rotation of planar bodies around a coplanar 
axis was studied also in \cite{DMV13}, where however only the case 
of an external axis (standard cylindrical coordinates) was considered. 
By using generalized cylindrical coordinates the rotation axis can 
intersect the interior of the body (a planar disk section here). 

Now, if the rotated planar domain, say $D$, is a blending domain (see 
Section 
2.1), then $r=r(t,\theta)$, $z=z(t,\theta)$ have the form 
(\ref{blending}), so that $r,z\in 
\mathbb{P}_1([0,1])\bigotimes \mathbb{T}_1([\omega_1,\omega_2])$, 
and the transformation $\sigma$ takes on the form
\begin{equation} \label{sigmarot}
\sigma(t,\theta,\phi)=(r(t,\theta)\cos(\phi),r(t,\theta)\sin(\phi),
z(t,\theta))\;,
\end{equation}
where $(t,\theta,\phi)\in [0,1]\times [\omega_1,\omega_2]\times 
[\phi_1,\phi_2]$, which clearly falls in the class of Proposition 1. 

This allows to construct optimal polynomial meshes on several 
arcwise solid sections, for example {\em solid caps} (one of the two 
portions of a ball  cut by a plane), and  
{\em spherical zones}  
(the portion of a ball between two parallel cutting planes), 
that correspond to a complete
rotation (i.e., by a multiple of $\pi$) of a planar circular segment or 
zone, respectively, around their symmetry axis. 

A rotation of a half-disk around the diameter 
by an angle smaller than $2\pi$  
produces a {\em spherical slice} (whose external boundary is a spherical 
lune).  

A {\em spherical cone} corresponds to a complete rotation of a 
planar circular sector around its axis, whereas a {spherical lens} (the 
intersection of 
two overlapping balls) to a complete 
rotation of a planar lens (around the axis connecting the centers), 
and a {\em spherical shell} to a complete rotation of a disk annulus 
around a diameter. 

Similarly, a complete rotation of circular segments and zones around an 
external axis 
parallel to their symmetry axis  
produce {\em solid toroidal caps} and {\em zones}, respectively 
(whose external boundaries are surface toroidal caps and collars).  

If the entire disk is rotated around an external axis by an angle smaller  
than $2\pi$, we get a {\em solid toroidal slice} (whose external boundary 
is a surface 
toroidal slice).  

In all the cases above the corresponding polynomial meshes
\begin{equation} \label{cylmesh}
\mathcal{A}_n(m)=\sigma\left(X_{mn}([0,1])\times 
W_{mn}([\omega_1,\omega_2])\times
W_{mn}([\phi_1,\phi_2])\right)\;,
\end{equation}
are optimal, with $c=\alpha^2\beta$, except for the 
spherical shell where we can 
take $[\omega_1,\omega_2]=[\phi_1,\phi_2]=[0,2\pi]$ and 
thus $c=\alpha^3$. In fact, $d_1+d_2+d_3=3$
and $N=dim(\mathbb{P}_n^3(K))=(n+1)(n+2)(n+3)/6\sim n^3/6$.  

If the rotated planar domain $D$ is a lune, then $r=r(\theta_1,\theta_2)$ 
and 
$z=z(\theta_1,\theta_2)$ have the form (\ref{lune}), 
so that $r,z\in
\mathbb{T}_1([\omega_1,\omega_2])\bigotimes 
\mathbb{T}_1([\omega_3,\omega_4])$, and the 
transformation 
becomes 
\begin{equation} \label{sigmasolidlune}
\sigma(\theta_1,\theta_2,\phi)=(r(\theta_1,\theta_2)\cos(\phi),
r(\theta_1,\theta_2)\sin(\phi),
z(\theta_1,\theta_2))\;,
\end{equation}
where $(\theta_1,\theta_2,\phi)\in [\omega_1,\omega_2]\times
[\omega_3,\omega_4]\times [\phi_1,\phi_2]$ are all angular variables. 
In case the rotation axis is the line connecting the two centers, 
we get a {\em solid lune} (the difference of a ball with a second 
overlapping ball). The optimal polynomial mesh is 
\begin{equation} \label{solidlunemesh}
\mathcal{A}_n(m)=\sigma\left(W_{mn}([\omega_1,\omega_2])\times
W_{mn}([\omega_3,\omega_4])\times
W_{mn}([\phi_1,\phi_2])\right)\;,
\end{equation}
with $c=\beta^3$.

A different situation, not of rotation type, arises when we consider 
a {\em spherical square pyramid}, that is a pyramid whose base is a 
geographic 
rectangle (the region of sphere comprised between two given latitudes
and longitudes) 
and whose vertex, 
say $V=(v_1,v_2,v_3)$, lies inside the ball. 
In this case, we can take the degenerate trivariate blending 
transformation 
\begin{equation} \label{sphersqpyr}
\sigma(t,\theta_1,\theta_2)=t\left(\cos(\theta_1)\cos(\theta_2),
\cos(\theta_1)
\sin(\theta_2),\sin(\theta_1)\right)+(1-t)V\;,
\end{equation}
where $(t,\theta_1,\theta_2)\in [0,1]\times [\omega_1,\omega_2]\times
[\omega_3,\omega_4]$. The optimal polynomial mesh has the form 
(\ref{cylmesh}) with $c=\alpha\beta^2$. 

In Table 1, we have listed (with no pretence of exhaustivity) the 
planar, surface and solid cases 
discussed above, together with the corresponding mesh parameters. 

\begin{table}
\begin{center}
\caption{Some standard planar, surface and solid sections 
together with  
the corresponding polynomial mesh 
parameters (cf. Proposition 1).}
\begin{tabular}{c c c}
$c$ & card. bound & section type\\
\hline 
$\alpha$ & $N_2$ & entire circle\\
\hline
$\beta$ & $N_2$ & circle arc\\
\hline
$\alpha^2$ & $N_1 N_2$ & entire disk/disk annulus\\
& $N_2^2$ & entire sphere/torus\\
\hline
$\alpha \beta$ & $N_1 N_2$ & disk sector/segment/zone/lens\\
& & surface spherical cap/collar\\
& & surface toroidal cap/collar/slice\\
\hline
$\beta^2$ & $N_2^2$ & planar lune\\
& & surface spherical rectangle/lune\\
& & surface toroidal rectangle\\
\hline
$\alpha^3$ & $N_1 N_2^2$ & entire ball/solid torus\\
& & spherical shell\\ 
\hline
$\alpha^2 \beta$ & $N_1 N_2^2$ & solid spherical cap/cone/lens/zone\\
& & solid toroidal cap/slice/zone\\
\hline
$\alpha \beta^2$ & $N_1 N_2^2$ & spherical square pyramid\\
\hline
$\beta^3$ & $N_2^3$ 
& solid lune\\
\hline
\end{tabular}
\end{center}
\end{table}


\small{}
\end{document}